  \documentclass[]{amsart}
   \newtheorem{lemma}{Lemma}[section]
   \newtheorem{theorem}[lemma]{Theorem}

   \newtheorem{coro}[lemma]{Corollary}
   \newtheorem{definition}[lemma]{Definition}

   \newcommand{\Wsob}{\smash{{\stackrel{\circ}{W}}}_2^1(D)}

\newcommand{\D}{\Delta}
\newcommand{\p}{\partial}
\renewcommand{\phi}{\varphi}

\renewcommand{\k}{\kappa}

\title[Enstrophy and Ergodicity of Gravity Currents]
{Enstrophy and Ergodicity of Gravity Currents }

\author{Vena Pearl Bongolan-Walsh, Jinqiao Duan}
\address[V. P. Bongolan-Walsh and J. Duan]
{Department of Applied Mathematics\\
 Illinois Institute of Technology\\
   Chicago, IL 60616, USA }
\email[]{bongven@iit.edu; duan@iit.edu}

\author{Hongjun Gao}
\address[H. Gao]
{Department of Mathematics\\
Nanjing Normal University\\
Nanjing 210097, China}
\email[H.~Gao]{gaohj@njnu.edu.cn}

\author{Tamay \"{O}zg\"{o}kmen }
\address[T. \"{O}zg\"{o}kmen]
{RSMAS/MPO\\University of Miami\\Miami, Florida, USA} \email[T.
 \"{O}zg\"{o}kmen]{tamay@rsmas.miami.edu}

\author{Paul Fischer}
\address[P. Fischer]{Argonne National Laboratory\\Argonne,
Illinois, USA} \email[P. Fischer]{fischer@mcs.anl.gov}

\author{Traian Iliescu}
\address[T. Iliescu]{Mathematics Department\\Virginia Tech\\
Blacksburg, VA 24061, USA}
\email[T.Iliescu]{iliescu@calvin.math.vt.edu}

\date{May 4, 2004}


\begin{document}

\begin{abstract}

    We study a coupled deterministic   system of vorticity evolution and salinity
    transport equations, with spatially correlated white noise on the boundary. This
    system may be considered as a model for gravity currents in
    oceanic fluids. The noise is due to   uncertainty in
    salinity flux on fluid boundary. After transforming this system into a random
    dynamical system, we first obtain an asymptotic estimate of enstrophy evolution,
    and then show that the system
    is ergodic under suitable conditions on mean salinity input flux on the boundary,
    Prandtl number and covariance of the noise.

\end{abstract}

\maketitle

\section{Geophysical background}

A gravity current is the flow of one fluid within another driven
by the gravitational force acting on the density difference
between the fluids. Gravity currents occur in a wide variety of
geophysical fluids.

Oceanic gravity currents are of particular importance, as they are
intimately related to the ocean's role in climate dynamics. The
thermohaline circulation in the ocean is strongly influenced by
dense-water formation that takes place mainly in polar seas by
cooling   and in marginal seas by evaporation. Such dense water
masses are released into the large-scale ocean circulation in the
form of overflows, which are bottom gravity currents.

 We   consider a two-dimensional
    model for oceanic gravity currents, in terms of   the
Navier-Stokes equation in vorticity form  and the transport
equation for salinity. The Neumann boundary conditions for this
model involve a spatially correlated white noise due to uncertain
salinity flux at the inlet boundary of the gravity currents.

In the  next section, we present the    model and reformulate it
  as a random dynamical system,  and then discuss the
cocycle property and dissipativity of this     model in \S 3 and
\S 4, respectively. Main results on random attractors, enstrophy
and ergodicity are in \S 5.  Enstrophy is   one half of the
mean-square spatial integral of vorticity. Ergodicity implies that
the time average for observables of the dynamical system
approximates the statistical ensemble average, as long as the time
interval is sufficiently long.

\section{Mathematical model}\label{sect2}

Oceanic gravity currents are usually down a slope of small angle
 (order of a few degrees). We model the gravity currents in the
  downstream-vertical plane, and ignore the variability in the
 cross-stream direction. This is an   appropriate approximate
 model for, e.g., the Red Sea overflow that flows along a long
 narrow channel that naturally restricts motion in the lateral
 planar plane \cite{Oz}. In fact, we will ignore small slope angle
 and the rotation, both affect the following estimates
 non-essentially, i.e., our results below still hold with
 non-essential modification of constants in the estimates and in the
 conditions for the ergodicity.

 Thus we   consider the gravity currents in the downstream-vertical
  $(x,z)$-plane. It is composed of
the Boussinesq equations for ocean fluid dynamics in terms of
vorticity $q(x,z,t)$, and the   transport  equation for
  oceanic salinity $S(x,z, t)$   on the domain $D=\{ (x,z):
0\leq x, z \leq 1 \}$:
\begin{align}\label{eqn1}
\begin{split}
q_t + J(q, \psi )  = &  \D q  -{ Ra}\partial_xS, \\
S_t + J(S, \psi )  = & {\rm 1 \over Pr}\D S,
\end{split}
\end{align}
where
\[
q(x,z, t)= -\D \psi
\]
is the vorticity in terms of stream function $\psi$,   $Pr$ is the
Prandtl number and Ra is the Rayleigh number. Moreover,
$J(g,h)=g_xh_z-g_zh_x$ is the Jacobian operator and
$\D=\p_{xx}+\p_{zz}$ is the Laplacian operator.   All these
equations are in non-dimensionalized forms. For the simplicity, we
let $Pr = 1$.

Note that the Laplacian operator $\D$ in the temperature and
salinity transport equations  is presumably $\p_{xx}+
\frac{\k_V}{\k_H} \delta^2 \p_{zz}$ with  $\delta$ being the
aspect ratio,  and $\k_H, \k_V$ the horizontal and vertical
diffusivities of salt, respectively.  However, our energy-type
estimates and the results below will not be essentially affected
by taking a homogenized  Laplacian operator $\D=\p_{xx}+\p_{zz}$.
All  our results would be true for this modified Laplacian. The
effect of the rotation   is parameterized in the magnitude of the
viscosity and diffusivity terms as discussed in \cite{ThuMcW92}.

The  fluid boundary condition is no normal flow and free-slip
on the whole boundary
\begin{equation*}
\psi=0,  \;  q = 0.
\end{equation*}

The flux boundary conditions are assumed for the ocean   salinity
$S$.

 At the inlet boundary $ \{x=0, \; 0 <z<1\}$ the flux is specified as:
\begin{equation}
\label{btop}
   \partial_x S =F(z)+ \dot{w}(z,t),
 \end{equation}
with $F(z)$ being the mean freshwater flux, and the  fluctuating
  part $\dot{w}(z, t)$   is
usually  of a shorter time scale than the response time scale of
the   oceanic  mean salinity. So we neglect the autocorrelation
time of this fluctuating process and thus assume that the noise is
white in time. The spatially correlated white-in-time noise
$\dot{w}(z, t)$ is described as the generalized time derivative of
a Wiener process $w(z, t)$ defined in a probability space
$(\Omega, \mathbb{F}, \mathbb{P})$, with mean vector zero and
covariance operator $Q$.

On the   outlet boundary $\{x=1,\; 0<z<1 \}$:
\begin{equation*}
\label{side}  \partial_x S =0.
\end{equation*}

At the top boundary $z=1$, and at the bottom boundary $z=0$:
\begin{equation*}
  \partial_z S=0  \;.
\end{equation*}

This is a system of deterministic   partial differential equations
with a stochastic boundary condition.

\section{Cocycle property }

In this section we will show that (\ref{eqn1})   has a unique
solution, and by   reformulating the model, we see it defines a
 cocycle or a random dynamical
system.

\smallskip

For the following we need some tools from the theory of partial differential equations.\\
Let $W_2^1(D)$ be the Sobolev space of functions on $D$ with first
generalized derivative in $L_2(D)$, the function space of square
integrable functions on $D$ with norm and  inner product
\[
\|u\|_{L_2}=\left(\int_D|u(x)|^2dD\right)^\frac{1}{2},\quad
(u,v)_{L_2}=\int_Du(x)v(x)dD,\quad u,\,v\in L_2(D).
\]
The space $W_2^1(D)$ is equipped with the norm
\[
\|u\|_{W_2^1}=\|u\|_{L_2}+\|\partial_x u\|_{L_2}+\|\partial_z
u\|_{L_2}.
\]
Motivated by the zero-boundary conditions of $q$ we also introduce the space
$\Wsob$ which contains roughly speaking  functions which are zero on the boundary $\partial D$
of $D$.
This space can be equipped with the norm
\begin{equation}\label{eqno}
\|u\|_{\Wsob}=\|\partial_x u\|_{L_2}+\|\partial_z u\|_{L_2}.
\end{equation}

Similarly, we can define function spaces on the interval $(0,1)$
denoted by $L_2(0,1)$.\\
Another Sobolev space is given by $\dot{W}_2^1(D)$ which is a
subspace of $W_2^1(D)$ consisting of functions $h$ such that
$\int_D hdD=0$. A norm equivalent to the $W_2^1$-norm on
$\dot{W}_2^1(D)$ is given by the right hand side of (\ref{eqno}).
For the subspace of functions in $L_2(D)$
having this property we will denote as $\dot{L}_2(D)$.\\


 We reformulate the above stochastic initial-boundary value problem
 into a random dynamical system \cite{Arn98}.  For convenience, we introduce
  vector notation for  unknown geophysical quantities
\begin{equation}\label{eqquant}
u=(q,S).
\end{equation}

\bigskip

Let $\dot{w}$ be a white noise in $L_2(0,1)$ with finite trace of
    the covariance operator  $Q$, and    the Wiener process $w(t)$
be defined on a probability space $(\Omega, \mathbb{F},
\mathbb{P})$.

\smallskip

Now we choose an appropriate phase space $H$ for this system. We
assume that the mean salinity flux $F\in L_2(0,1)$. Note that
\[\frac{d}{dt}\int_D S dD = \int_0^1 [F(z)+ \dot{w}(z,t)] dz = \mbox{constant}.
\]
It is reasonable (see \cite{Dij00})  to assume that
\begin{equation}
\int_0^1 [F(z)+ \dot{w}(z,t)] dz = 0,
\end{equation}
and thus $ \int_DSdD$ is constant in time and, without loss of
generality, we may assume it is zero (otherwise, we subtract this
constant from $S$):
\begin{equation*}
 \int_DSdD = 0.
\end{equation*}
Thus $S\in \dot{L}_2(D)$, and we have the usual Poincar\'e
inequality for $S$.

\bigskip
Define the phase space
$$
H=L_2(D) \times \dot{L}_2(D).
$$
We rewrite the coupled   system (\ref{eqn1}) as:
\begin{equation} \label{eqn19}
\frac{du}{dt}+Au=F_1(u) + F_2(u),\qquad u(0)=u_0\in H,
\end{equation}

\smallskip

where
\begin{equation*}
Au =\left(
\begin{array}{l}
-\Delta q \\
-{\rm 1 \over Pr}\Delta S
\end{array}
\right),
\end{equation*}

\smallskip

\[
F_1(u)[x,z]= \left(
\begin{array}{c}
-J(q,\psi)\\
-J(S,\psi)
\end{array}
\right)[x,z],
\]
and
\[
F_2(u)[x,z]= \left(
\begin{array}{c}
\,{\rm Ra}(-\partial_x S)\\
0
\end{array}
\right)[x,z].
\]

The boundary and initial conditions   are:
\begin{align}\label{eqn199}
\begin{split}
q  = & 0, \mbox{on}\; \partial D,\;  \\
\frac{\partial S}{\partial n}  = & 0, \;\mbox{on}\; \partial
D\backslash\{x = 0,\; 0 <z<1\}, \; \\
\frac{\partial S}{\partial x}  = & F(z)+
\dot{w}(z,t),\;\mbox{on}\;
\{x =0, \;0 <z<1\},\; \\
u(0)=& u_0= \left(
\begin{array}{c}
q_o\\
S_o
\end{array}
\right),
\end{split}
\end{align}
where $n$ is the out unit normal vector of $\partial D$. The
system (\ref{eqn19}) consists of deterministic partial
differential equations with stochastic Neumann boundary
conditions.
\smallskip

\bigskip

We now transform the above system (\ref{eqn19}) into a system of
random partial differential equations (i.e., evolution equations
with random coefficients) with homogeneous boundary conditions,
whose solution map can be easily seen as a cocycle. Thus we can
investigate this dynamics in the framework of random dynamical
systems \cite{Arn98}. Note that we have a non-homogenous
stochastic boundary condition for salinity $S$, so the first step
is to homogenize this boundary condition.



\bigskip

To this end, we need  an Ornstein-Uhlenbeck stochastic process
solving the linear  differential equation

\begin{equation}\label{eq20}
\frac{d\eta_1}{dt}  =  \Delta \eta_1,
\end{equation}
with the following same boundary conditions as for $S$, and zero
initial condition
\begin{align}\label{eq200}
\begin{split}
 \partial_x \eta_1(t,0,z,\omega)=& F(z)+ \dot{w}(z,t), \\
 \partial_x \eta_1(t,1,z,\omega)=& 0, \\
 \partial_z \eta_1(t,x,0,\omega)=& 0, \\
 \partial_z \eta_1(t,x,1,\omega)=& 0, \\
 \eta_1(0,x,z,\omega) = & 0, \\
 \int_D \eta_1 dD = & 0.
\end{split}
\end{align}

\smallskip

\begin{lemma}
Suppose that the covariance $Q$ has  finite trace : ${\rm
tr}_{L_2} Q<\infty$. Then the Ornstein-Uhlenbeck problem
(\ref{eq20})-(\ref{eq200})   has a unique stationary solution in
$L_2(D)$ generated by
\[
(t,\omega)\to \eta_1(\theta_t\omega).
\]

\end{lemma}

In fact, we can write down the solution $\eta_1$ following Da
Prato and Zabczyck \cite{DaPZab96, DPZ93} as

\begin{equation}
\label{btop}
 \eta_1(t,x,z,\omega) =( - \Delta) \int_0^t S(t-s) {\mathcal{N} (\textbf{X}) } ds
 \end{equation}

\smallskip

    where $I$ is the identity operator in $L_2(D)$, and
    $\mathcal{N}$ is the solution operator to the elliptic
    boundary value
    problem $ \Delta h -\lambda h= 0$ with the boundary conditions for $h$     the
     same as $\eta_1$, that is
    $\partial h/\partial n = \textbf{X}$ on $\partial D$ with
    $\int_D h dD = 0$, where $n$ is the unit outer normal vector to
    $\partial D$ and $\textbf{X}$ is
\[
\textbf{X} = \left(
\begin{array}{c}
F(z)+ \dot{w}(z,s)\\
0\\
0\\
0\\
\end{array}
\right).
\]
Here $\lambda$ is chosen so that this elliptic
    boundary value  problem has a unique solution. Since $\int_D h dD = 0$,
we can choose $\lambda = 0$. Moreover, S(t) is a strongly
continuous semigroup, symbolically, $e^{\Delta t}$, that is, the
generator of $S(t)$ is $\Delta$.

\bigskip

Now we are ready to transform (\ref{eqn19}) into a random
dynamical system in Hilbert space $H$. Define
\begin{equation}\label{eq22}
\eta(t,x,z,\omega)= \left(
\begin{array}{c}
0\\
\eta_1
\end{array}
\right)
\end{equation}
and recall
\begin{equation*}
u =\left(
\begin{array}{l}
q\\
S
\end{array}
\right).
\end{equation*}

Let

\begin{equation}\label{eq19a}
v:=u - \eta.
\end{equation}

Then we obtain a random partial differential equation

\begin{equation}\label{eq18}
\frac{dv}{dt}+Av=F_1(v+\eta(\theta_t\omega))+F_2(v+\eta(\theta_t\omega)),\qquad
v(0)=v_0\in H .
\end{equation}

\bigskip

where

\smallskip

$v_x(t,0,z,\omega)$, $v_x(t,1,z,\omega)$, $v_z(t,x,0,\omega)$ and
$v_z(t,x,1,\omega)$ are now all zero vectors,i.e., homogenous
boundary conditions, and the  initial condition is still the same
as for $u$

\smallskip

\[
v(0,x,z,\omega)= \left(
\begin{array}{c}
q_o\\
S_o
\end{array}
\right),
\]

\smallskip


\smallskip

However, because of the Jacobian, we will have to do nonlinear
analysis on (\ref{eq18}) to resolve v.

\bigskip


We  introduce another   space

\smallskip

$$V=\Wsob\times  \dot{W}_2^1(D)$$.

\smallskip

For sufficiently smooth functions $v=(\bar q,\bar S)$, we can
calculate via integration by parts

\begin{align}\label{eq1}
\begin{split}
({\mathcal{A}}u,v)_{H} = &
\int_D\nabla q\cdot\nabla \bar q\,dD  \\
& + {1 \over Pr}\int_D\nabla S\cdot\nabla \bar S\,dD\,
\end{split}
\end{align}

\smallskip

since now, we have only homogenous boundary conditions.

\smallskip

 Hence on the  space $V$ we can introduce a bilinear
form $\tilde{a}(\cdot, \cdot)$ which is continuous, symmetric and
positive
\begin{align*}
\tilde{a}(u,v)=& \int_D\nabla q\cdot\nabla \bar qdD   +
\int_D\nabla S\cdot\nabla \bar SdD
\end{align*}
This bilinear form defines a unique linear continuous operator
$A:V\to V^\prime$ such that $\langle Au,v\rangle=\tilde{a}(u,v)$.

\smallskip

Recall

\smallskip

\[
F_1(u)[x,z]= \left(
\begin{array}{c}
-J(q,\psi)\\
-J(S,\psi)
\end{array}
\right)[x,z].
\]
and
\[
F_2(u)[x,z]= \left(
\begin{array}{c}
\,{\rm Ra}(-\partial_x S)\\
0
\end{array}
\right)[x,z].
\]

\bigskip

\begin{lemma}\label{lcon}
The operator $F_1:V\to H$ is continuous. In particular, we have
\[
\langle F_1(u),u\rangle=0.
\]
\end{lemma}
\begin{proof}
We have a constant $c_1>0$ such that
\begin{equation}\label{eq6a}
\|\psi\|_{W_2^3(D)}\le c_1\|q\|_{W_2^1(D)}
\end{equation}
for any $q\in W_2^1(D)$ which follows straight forwardly by
regularity properties of a linear elliptic boundary problem. Note
that $W_2^3$ is a Sobolev space with respect to the third
derivatives. Hence we get:
\begin{equation*}
\begin{split}
\|J(S,\psi)\|_{L_2} &
\le \sup_{(x,z)\in D}(|\partial_x\psi(x,z)|+|\partial_z\psi(x,z)|)\times\\
&
\times\left(\int_D|\partial_xS(x,z)|+|\partial_zS(x,z)|dD\right).
\end{split}
\end{equation*}
  The second factor on the right hand side is bounded by

\smallskip

\[
\left(\int_D|\partial_xS(x,z)|^2dD\right)^\frac{1}{2}
+\left(\int_D|\partial_zS(x,z)|^2dD\right)^\frac{1}{2}\le \|u\|_V.
\]

\smallskip

On account of the Sobolev embedding Lemma,  we have some positive
constants  $c_2,\, c_3$ such that
\[
\sup_{(x,z)\in D}(|\partial_x\psi(x,z)|+|\partial_z\psi(x,z)|) \le
c_2 \|\nabla \psi\|_{W_2^2(D)}\le c_3\|q\|_{W_2^1(D)}\le c_3
\|u\|_V.
\]
Hence we have a positive constant $c_4$ such that
\[
\|J(S,\psi)\|_{L_2}\le c_4\|u\|_V^2
\]
for $u\in V$.

\smallskip

We now show that
\[
\langle J(S,\psi),S\rangle =0.
\]

We obtain via integration by parts
\begin{align*}
\int_D\partial_xS&\partial_z\psi\, SdD-\int_D\partial_zS\partial_x\psi\, SdD\\
=&-\int_D\partial_{xz}^2S\psi \,SdD+\int_D\partial_{zx}^2S\psi
S\,dD
-\int_D\partial_x S\psi\partial_zS\,dD+\int_D\partial_z S\psi\partial_xS\,dD\\
&+\int_{(0,1)}\partial_xS\psi S|_{z=0}^{z=1}dx-
\int_{(0,1)}\partial_zS\psi S|_{x=0}^{x=1}dz=0
\end{align*}
because $\psi$ is zero on the boundary $\partial D$. This relation
is true for a set of  sufficiently smooth functions $\psi,\,S$
which are dense in $\Wsob\times W_2^1(D)$. By the continuity of
$F_1$, as just shown in Lemma \ref{lcon}, we can extend this
property to $\Wsob\times W_2^1(D)$.
\end{proof}

\bigskip

\begin{lemma}
The following estimate holds
\[
\|F_2(u)\|_{L^2}\le c_5\|u\|_V.
\]
 for some  positive constant $c_5$.
\end{lemma}

\smallskip

\begin{proof}
By simple calculation, the proof is obtained.
\end{proof}



\smallskip

We have obtained a differential equation without white noise but
with random coefficients. Such a differential equation can be
treated sample-wise for {\em any} sample $\omega$. We are looking
for solutions in
\begin{equation*}
v\in C([0, \tau]; H)\cap L^2(0, \tau; V),
\end{equation*}
for all $\tau>0$. If we can solve this equation then $u:=v+\eta$
defines a solution version of (\ref{eqn19}).

For the well posedness of the problem we now have the following
result.

\begin{theorem}\label{tEX}
({\bf Well-Posedness})
For any time $\tau>0$, there exists a unique solution of (\ref{eq18})
in $C([0,\tau];H)\cap L_2(0,\tau;V)$. In particular, the solution mapping
\[
{{\mathbb{R}}}^+\times \Omega\times H\ni(t,\omega,v_0)\to v(t)\in H
\]

is measurable in its arguments and the solution mapping
$H\ni v_0\to v(t)\in H$ is continuous.\\
\end{theorem}
\begin{proof}
By the properties of $A$ and $F_1$ (see Lemma \ref{lcon}),
the random   differential equation  (\ref{eq18})  is   essentially similar  to the
2 dimensional  Navier Stokes equation.  Note that $F_2$ is only an affine mapping.
Hence we have existence and uniqueness
and the above regularity assertions.
\end{proof}

On account of the transformation (\ref{eq19a}),  we find that
(\ref{eqn19}) also has a unique solution.

Since the solution mapping
\[
{\mathbb{R}}^+\times\Omega\times H\ni (t,\omega,v_0)\to v(t,\omega,v_0)=:\phi(t,\omega,v_0)\in H
\]
is well defined,  we can introduce a random dynamical system. On $\Omega$ we can define
a shift operator $\theta_t$ on the paths of the Wiener process
that pushes our noise:
\[
w(\cdot,\theta_t\omega)=w(\cdot+t,\omega)-w(t,\omega)\quad \text{for }t\in{\mathbb{R}}
\]
which is called the {\em Wiener shift}.
Then $\{\theta_t\}_{t\in{\mathbb{R}}}$ forms a flow which is ergodic for the probability measure
${\mathbb{P}}$.
The properties of the solution mapping cause  the following relations
\begin{align*}
&\phi(t+\tau,\omega,u)=\phi(t,\theta_\tau\omega,\phi(\tau,\omega,u))\quad\text{for }t,\,\tau\ge 0\\
&\phi(0,\omega,u)=u
\end{align*}
for any $\omega\in\Omega$ and $u\in H$. This property is called the cocycle property of $\phi$
which is important to study the dynamics of random systems. It is a generalization
of the semigroup property. The cocycle $\phi$ together with the flow $\theta$ forms
a {\em random dynamical system}.

\section{Dissipativity}

In this section we   show that   the  random dynamical system
(\ref{eq18}) for gravity currents
 is dissipative,  in the sense that it has an absorbing  (random)
 set.
This means that the solution $v$ is contained in a particular
region of the phase space $H$ after a sufficiently long time.
 This dissipativity will help us to obtain
asymptotic estimates of the  enstrophy and salinity evolution.
Dynamical properties that follow  from this dissipativity will be
considered in the next section. In particular, we will show that
the   system has a random attractor, and is ergodic
under suitable conditions.\\

\smallskip

We introduce the spaces
\begin{align*}
\tilde H&= \dot{L}_2(D)\\
\tilde V&=\dot{W}_2^1(D).
\end{align*}
We also choose a subset of  dynamical variables of our system (\ref{eqn1}).
\begin{equation}
\label{eqE} \tilde v=   S-\eta_1 .
\end{equation}
 To calculate the energy inequality
for $\tilde v$,   we apply the chain rule to $\|\tilde v\|_H^2$.
We obtain by Lemma \ref{lcon}
\begin{align}\label{eqE1}
\begin{split}
\frac{d}{dt} &\|\tilde v\|_{\tilde H}^2+2\|\nabla\tilde v\|_{L_2}^2\\
=&2\langle J(\eta_1, \psi), \tilde v\rangle.
\end{split}
\end{align}
 The
expression $\nabla \tilde v$ is defined by $(\nabla_{x,z}\tilde
v)$. We now can estimate  the term on the right hand side.

\smallskip

 By the Cauchy inequality,  integration by parts and Poincare inequality $\lambda_1\|q\|_{L_2}\le \|\nabla q\|_{L_2}$ for $q\in
\Wsob$ and $\lambda_2\|\tilde v\|_{L_2}\le \|\nabla \tilde
v\|_{L_2}$ for $\tilde v\in \dot{W}_2^1(D)$, we have
 \begin{equation}\label{eqv}
\frac{d}{dt} \|\tilde v\|_{\tilde H}^2 + \|\nabla\tilde
v\|_{L_2}^2 \le \lambda_1^2\|\nabla q\|^2\|\eta_1\|^2.
\end{equation}
For $q$, we have the following estimate
\begin{equation}\label{eqq}
\frac{d}{dt}\|q\|^2 + \|\nabla q\|^2 \le
Ra^2(\lambda_2^2\|\nabla\tilde v\|^2 + \|\eta_1\|^2).
\end{equation}
From (\ref{eqv}) and (\ref{eqq}), we have
\begin{equation}\label{eqvq}
\frac{d}{dt} (2\|\tilde v\|_{\tilde H}^2 +
\frac{1}{Ra^2\lambda_2^2}\|q\|^2) + \|\nabla\tilde v\|_{L_2}^2 +
(\frac{1}{Ra^2\lambda_2^2} - 2\lambda_1^2\|\eta_1|^2)\|\nabla
q\|^2\le \frac{1}{\lambda_2^2}\|\eta_1\|^2.
\end{equation}







\begin{definition}\label{defA}
A random set $B=\{B(\omega)\}_{\omega\in \Omega}$ consisting of closed bounded sets $B(\omega)$
is called absorbing for a random dynamical system $\phi$ if we have for
any  random set $D=\{D(\omega)\}_{\omega\in\Omega},\,D(\omega)\in H$ bounded,
 such that $t\to \sup_{y\in D(\theta_t\omega)}\|y\|_H$ has a subexponential
growth for $t\to\pm\infty$
\begin{align}\label{eqA0}
\begin{split}
&\phi(t,\omega,D(\omega))\subset B(\theta_t\omega) \quad\text{for }t\ge t_0(D,\omega)\\
&\phi(t,\theta_{-t}\omega,D(\theta_{-t}\omega))\subset B(\omega)\quad\text{for }t\ge t_0(D,\omega).
\end{split}
\end{align}
$B$ is called forward invariant if
\[
\phi(t,\omega,u_0)\in B(\theta_t\omega)\quad \text{if }  u_0\in B(\omega)\quad \text{for }t\ge 0.
\]
\end{definition}
Although $\tilde v$ is not a random dynamical system in the strong sense
we can also show dissipativity in the sense of the above definition.

\begin{lemma}\label{lemA}
Let $\phi(t,\omega,v_0)\in H$ for $v_0\in H$ be defined in
(\ref{eqn19}), and
 $$
 \frac{1}{Ra^2\lambda_2^2} - 2\lambda_1^2{\mathbb{E}}\|\eta_1|^2 > 0.
 $$
 Then the closed ball $B(0,R_1(\omega))$ with
radius
\[
R_1(\omega)=
2\int_{-\infty}^0e^{\alpha\tau}\frac{1}{\lambda_2^2}\|\eta_1\|^2
d\tau
\]
is forward invariant and absorbing.
\end{lemma}

The proof of this lemma follows by integration of (\ref{eqvq}).\\

For the applications in the next section we need that the elements which are contained in the
absorbing set satisfy a particular regularity. To this end we introduce the function space
\[
{\mathcal{H}}^s:=\{u\in H: \|u\|_s^2:=\|A^\frac{s}{2}u\|_H^2<\infty\}
\]
where $s\in {\mathbb{R}}$. The operator $A^s$ is the $s$-th power of the positive
and symmetric operator $A$.
Note that these spaces are embedded in the Slobodeckij spaces
$H^s,\, s>0$. The norm of these spaces is denoted by $\|\cdot\|_{H^s}$.
This norm can be found in Egorov and Shubin \cite{EgoShu91}, Page 118.
But we do not need this norm explicitly. We only mention
that on ${\mathcal{H}}^s$ the norm
$\|\cdot\|_s$
of
$H^s$ is equivalent to the norm of ${\mathcal{H}}^s$ for $0 < s$, see \cite{LioMag68}.\\


Our goal is it to show that  $v(1,\omega,D)$ is  a bounded set in
${\mathcal{H}}^s$ for some $s>0$.
This property causes the complete continuity of the mapping $v(1,\omega,\cdot)$.
We now derive a differential inequality for $t\|v(t)\|_s^2$.
By the chain rule we have
\[
\frac{d}{dt}(t\|v(t)\|_s^2)=\|v(t)\|_s^2+t\frac{d}{dt}\|v(t)\|_s^2.
\]
Note that for the embedding  constant $c_{s}$ between
${\mathcal{H}}^s$ and $V$
\[
\int_0^t\|v\|_s^2ds\le c_{s}^2 \int_0^t\|v\|_V^2ds\quad \text{for
} s \le 1
\]
such that the left hand  side is bounded if the initial conditions $v_0$
are contained in a bounded set in $H$.
The second term in the above formula can be expressed as followed:
\begin{align*}
t\frac{d}{dt}(A^\frac{s}{2}v,A^\frac{s}{2}v)_H=&2t(\frac{d}{dt}v,A^sv)_H
=-2t(Au,A^sv)_H+2t(F_1(v+\eta(\theta_t\omega)),A^sv)_H\\
&+2t(F_2(v+\eta(\theta_t\omega)),A^sv)_H.
\end{align*}
We have
\[
(Av,A^sv)_H=\|A^{\frac{1}{2}+\frac{s}{2}}v\|_H=\|v\|_{1+s}^2.
\]
Similar to the argument of \cite{duansch} and the estimate for the
existence of absorbing, and applying some embedding theorems, see
Temam \cite{Tem83} Page 12 we have got
\begin{align*}
\|v\|_s^2 \le C(t, \|v_0\|_H, \sup\limits_{t\in [0,
1]}\|\eta_1\|_{D(A^s)}),\; \mbox{for}\; t\in [0, 1].
\end{align*}
By the results of \cite{DPZ93} and \cite{duansch}, we know $$
\sup\limits_{t\in [0, 1]}\|\eta_1\|_{D(A^s)} \le C(tr_{L^2}Q) <
\infty,\; \mbox{for}\; 0 < s < \frac14.
$$




By now, our estimates allow us to write down the main assertion
with respect to the dissipativity of this section.
\begin{theorem}\label{tA}
For the random dynamical system generated by (\ref{eq18}),
there exists a compact random  set $B=\{B(\omega)\}_{\omega\in \Omega}$
which satisfies Definition \ref{defA}.
\end{theorem}
We define
\begin{equation}\label{eqAB}
B(\omega)=\overline{\phi(1,\theta_{-1}\omega,B(0,R(\theta_{-1}\omega)))}\subset
{\mathcal{H}}^s,\quad 0<s<\frac14.
\end{equation}
In particular, ${\mathcal{H}}^s$ is compactly embedded in $H$.

\section{Random dynamics: Enstrophy and ergodicity}

In this section we will apply the  dissipativity result of the
last section to analyse the dynamical behavior of  the random
dynamical system (\ref{eqn1}). However, it will be enough to
analyse the transformed  random dynamical system generated by
(\ref{eq18}). By the transformation (\ref{eq19a}), we have
all these qualitative properties to the original gravity current system (\ref{eqn1}).\\

We will consider the following dynamical behavior: random
attractors, asymptotic evolution of enstrophy and mean-square norm
of salinity profile, and ergodicity. Enstrophy is defined as one
half of the mean-square integral of vorticity. Ergodicity implies
that the time average for observables of the dynamical system
approximates the statistical ensemble average, as long as the time
interval is sufficiently long.

We first consider random   attractors. We recall the following
basic concept; see, for instance, Flandoli and Schmalfu{\ss}
\cite{FlaSchm95a}.

\begin{definition}
Let $\phi$ be a random dynamical system. A random set
$A=\{A(\omega)\}_{\omega\in \Omega}$ consisting of compact
nonempty sets $A(\omega)$ is called random global attractor if for
any random bounded set $D$ we have for the limit in probability
\[
({\mathbb{P}})\lim_{t\to\infty} {\rm dist}_H(\phi(t,\omega,D(\omega)),A(\theta_t\omega)) = 0
\]
and
\[
\phi(t,\omega,A(\omega))=A(\theta_t\omega)
\]
any $t\ge 0$ and $\omega\in\Omega$.
\end{definition}

The essential long-time behavior of a random system is captured by
a random  attractor. In the last section we showed that the
dynamical system $\phi$ generated by (\ref{eqn19}) is dissipative
which means that there exists a random set $B$ satisfying
(\ref{eqA0}). In addition,  this set is compact. We now recall and
adapt the following theorem  from \cite{FlaSchm95a}.

\begin{theorem}
Let $\phi$ be a random dynamical dynamical system on the state
space $H$ which is a separable Banach space such that
$x\to\phi(t,\omega,x)$ is continuous. Suppose that $B$ is a set
ensuring the dissipativity given in Definition \ref{defA}. In
addition, $B$ has a subexponential growth (see  Definition
\ref{defA})  and is regular (compact). Then the dynamical system
$\phi$ has a random attractor.
\end{theorem}

This theorem can be applied to the random dynamical system
 $\phi$ generated by the stochastic differential equation (\ref{eq18}).
 Indeed, all the assumptions
are satisfied. The set $B$ is defined in Theorem \ref{tA}. Its
subexponential growth follows from $B(\omega)\subset
B(0,R(\omega))$ where the radius $R(\omega)$ has been introduced
in the last section. Note that $\phi$ is a {\em continuous} random
dynamical system; see Theorem \ref{tEX}. Thus $\phi$ has a random
attractor. By the transformation (\ref{eq19a}), this is also true
for the original gravity currents system (\ref{eqn19}).

\begin{coro}  \label{attractor}
({\bf Random Attractor})
 The  gravity currents system (\ref{eqn19})
has a random attractor.
\end{coro}

\bigskip

\bigskip

Now we   consider random fixed point (stationary state) and
ergodicity.

\begin{definition}
A random variable $v^\ast:\Omega\to H$ is defined to be a random fixed point for a random dynamical system
if
\[
\phi(t,\omega,v^\ast(\omega))=v^\ast(\theta_t\omega)
\]
for $t\ge 0$ and $\omega\in\Omega$. A random fixed point
$v^\ast$ is called exponentially attracting
if
\[
\lim_{t\to\infty}\|\phi(t,\omega,x)-v^\ast(\theta_t\omega)\|_H=0
\]
for any $x\in H$ and $\omega\in \Omega$.
\end{definition}
Sufficient conditions for the existence of random fixed points
are given in Schmalfu{\ss}
\cite{Schm97a}. We here formulate a simpler version of this
theorem and it is appropriate  for our system here.

\begin{theorem} ({\bf Random Fixed Point Theorem})
Let $\phi$ be a random dynamical system and suppose that $B$ is a forward invariant complete set.
In addition, $B$ has a subexponential growth, see Definition \ref{defA}.
Suppose that the following contraction conditions holds:
\begin{equation}\label{eqKX}
\sup_{v_1\not= v_2\in B(\omega)}\frac{\|\phi(1,\omega,v_1)-\phi(1,\omega,v_2)\|_H}{\|v_1-v_2\|_H}\le k(\omega)
\end{equation}
where the expectation of $\log k$ denoted by ${\mathbb{E}}\log k<0$.
Then $\phi$ has a unique random fixed point in $B$ which is exponentially attracting.
\end{theorem}

This theorem can be considered as a random version of the Banach fixed point theorem.
The contraction condition is formulated in the mean for the right hand side of  (\ref{eqKX}).\\

\begin{theorem} \label{fpt}
({\bf Unique Random Stationary State}) Assume  that the salinity
boundary flux data $\|F\|_{L_2}$, the Prandtl number $Pr$, and the
trace of the covariance for the noise ${\rm tr}_{L_2} Q$ are
sufficiently small. Then the random dynamical system generated by
(\ref{eqn19}) has a unique exponentially attracting random
stationary state.
\end{theorem}

Note that if we take into account the horizontal and vertical
diffusivities $\kappa_H, \kappa_V$ in the Laplacian operator, the
above ``smallness" condition needs to be slightly modified.

Here we only give a short sketch of the proof. Let us suppose for
a while that $B$ is given by the ball $B(0,R)$ introduced in Lemma
\ref{lEX}. Suppose that the data in the assumption of the lemma
are small and $\nu$ is large. Then it follows that ${\mathbb{E}}R$
is also small. To calculate the contraction condition we have to
calculate
 $\|\phi(1,\omega,v_1(\omega))-\phi(1,\omega,v_2(\omega))\|_H^2$
for arbitrary random variables  $v_1,\,v_2\in B$.
By Lemma \ref{lcon} we have that
\[
\langle J(q_1,\psi_1)-J(q_1,\psi_1),q_1-q_2\rangle\le
c_{22}\|q_1-q_2\|_{W_2^1}^2+c_{23}\|q_1\|_{\Wsob}^2\|q_1-q_2\|_{L_2}^2
\]
where the constant $c_{23}$ can be chosen sufficiently small if $\nu$ is large. On account
of the fact that the other expressions allow similar estimates and that $F_2$ is linear
we obtain:
\begin{align*}
\frac{d}{dt}&\|\phi(t,\omega,v_1(\omega))-\phi(t,\omega,v_2(\omega))\|_H^2\\
&\le
(-\alpha^\prime+c_{23}\|\phi(t,\omega,v_2(\omega))\|_V^2)\|\phi(t,\omega,v_1(\omega))-\phi(t,\omega,v_2(\omega))\|_H^2
\end{align*}
for some positive $\alpha^\prime $ depending on $\nu$.
From this inequality and the Gronwall lemma it follows that the contraction condition (\ref{eqKX}) is satisfied if
\[
{\mathbb{E}}\sup_{u_2\in
B(\omega)}c_{23}\int_0^1\|\phi(t,\omega,v_2)\|_V^2dt<\alpha^\prime.
\]
But by the energy inequality this property is satisfied if the ${\mathbb{E}}R$ and
${\mathbb{E}}\|z\|_V^2$ is sufficiently small which follows from the assumptions.
\\

Let now $B$ be the random set defined in (\ref{eqAB}).
Since the set $B$ introduced in (\ref{eqAB}) is absorbing any state the fixed point
$v^\ast$ is contained in this $B$. In addition $v^\ast$ attracts {\em any} state
from $H$ and not only states from $B$.

The uniqueness of  this random fixed point implies {\em ergodicity}.
We will comment on this issue at the end of this section.

\bigskip

By the well-posedness  Theorem \ref{tEX}, we know that the
stochastic evoltuion equation (\ref{eqn19})  has  a unique
solution. The solution is a Markov process.  We can define the
associated Markov operators $\mathcal T(t)$ for $t \geq 0$, as
discussed in \cite{Schm91a, Schm89}. Moreover, $\{{\mathcal
T}(t)\}_{t\geq 0}$ forms a semigroup.

Let $M^2$ be the set of probability distributions $\mu$ with finite energy,
i.e.,
\[
\int_H\|u\|_H^2d\mu(u)<\infty .
\]
Then the distribution of the solution $u(t)$ (at time $t$) of the
stochastic evolution equation   (\ref{eqn19}) is given by
\[
{\mathcal T}(t)\mu_0,\
\]
where the distribution $\mu_0$ of the initial data is contained in $M^2$.

We note that
the expectation of the solution $\|u(t)\|_H^2$ can be expressed in terms of this
distribution ${\mathcal T}(t)\mu_0$:
\[
{\mathbb{E}}\|u(t)\|_H^2=\int_H \|u\|_H^2d{\mathcal T}(t)\mu_0.
\]
We can derive the following energy inequality in the mean, using
our earlier estimates.

\begin{coro}\label{tMar}
The dynamical quantity  $u=( q,S)$ of the coupled
atmospheric-ocean system satisfy the estimate
\[
{\mathbb{E}}\|u(t)\|_H^2+\alpha{\mathbb{E}}\int_0^t\|u(\tau)\|_V^2d\tau
\le {\mathbb{E}}\|u_0\|_H^2+t\,c_{24}+t\,{\rm tr}_{L_2}Q,
\]
where the positive constants $c_{24}$ and $\alpha$ depend    on
physical data  $F(x)$,  $Pr$ and  $Ra$.
\end{coro}

By the Gronwall inequality, we further obtain the following result
about the asymptotic mean-square estimate. Recall that $\frac12
\mathbb{E} \int_D q^2dD$ is the {\em enstrophy} for gravity
currents. The following theorem gives asymptotic estimate for both
enstrophy and mean-square salinity $\mathbb{E} \int_D S^2dD$.

\begin{theorem} \label{feedbacktheorem}
 ({\bf Enstrophy and mean-square salinity})
For the expectation of the dynamical quantity  $u=( q, S)$ of the
coupled atmospheric-ocean system, we have the asymptotic estimate
 \[
\limsup_{t\to\infty}
{\mathbb{E}}\|u(t)\|_H^2=\limsup_{t\to\infty}\int_H\|u\|_H^2d{\mathcal  T}(t)\mu_0\le
\frac{c_{24}+{\rm tr}_{L_2}Q}{c_{25}}
\]
if the initial distribution $\mu_0$ of the random initial
condition $u_0(\omega)$ is contained in $M^2$.  Here   $c_{25}>0$
  depends only on physical data.

\end{theorem}

By the estimates of Corollary \ref{tMar},
 we are able to use the well known
Krylov-Bogolyubov procedure to conclude the existence of invariant
measures of the Markov semigroup $\mathcal T(t)$.

\begin{coro}
The semigroup of Markov operators $\{{\mathcal T}(t)\}_{t\geq 0}$ possesses an
invariant distribution $\mu_i$
in $M^2$:
\[
{\mathcal T}(t)\mu_i=\mu_i\quad \text{for }t\ge 0.
\]
\end{coro}

In fact,  the limit points of
\[
\left\{
\frac1t\int_0^t{\mathcal  T}(\tau)\mu_0d\tau\right\}_{t\geq 0}
\]
for $t\to\infty$ are invariant distributions. The existence of
such limit points follows from the estimate in Corollary
\ref{tMar}.

In some situations, the invariant measure may be unique. For
example, the  unique random fixed point  in  Theorem \ref{fpt} is
 defined by a random variable
$u^\ast(\omega)=v^\ast(\omega)+\eta(\omega)$. This random variable
corresponds to a unique invariant measure of the Markov semigroup
$\mathcal T(t)$. More specifically,
 this unique  invariant measure  is the expectation of the Dirac measure
with the random variable
as the random mass point
\[
\mu_i={\mathbb{E}}\delta_{u^\ast(\omega)}.
\]
Because the uniqueness of   invariant measure
implies   ergodicity \cite{DaPZab96}, we conclude that the
   gravity currents model (\ref{eqn1})
is ergodic under the   suitable conditions  in  Theorem \ref{fpt}
for physical data and random noise.  We reformulate it as the
following  ergodicity principle.

\begin{theorem}  \label{ergodic}
({\bf Ergodicity}) Assume that the salinity boundary flux data
$\|F\|_{L_2}$, the Prandtl number $Pr$, and the trace of the
covariance for the noise ${\rm tr}_{L_2}Q$ are sufficiently small.
Then the gravity currents system (\ref{eqn1}) is ergodic, namely,
for any observable of the gravity currents, its time average
approximates the statistical ensemble average, as long as the time
interval is sufficiently long.
\end{theorem}

In the regime of ergodicity, the gravity currents system can be
numerically simulated with (almost surely) one random sample. This
is what we call ergodicity-based numerical simulation of
stochastic systems, in contrast to sample-wise (``Monte Carlo")
simulations.

\bigskip

{\bf Acknowledgement.} This work was partly supported by the NSF
Grants DMS-0209326 and DMS-0139073,   and a Grant  of the NNSF of
China. H. Gao  would like to thank the Illinois Institute of
Technology, Chicago,
 for   hospitality.


\end{document}